\documentclass[a4paper,12pt]{article}
\usepackage{graphicx}
\usepackage{geometry}
\usepackage{amsmath}
\usepackage{amsfonts}
\usepackage{amssymb}
\usepackage[english]{babel}
\usepackage{multirow}
\usepackage{multicol}
\usepackage{amsthm}
\usepackage{subfigure}
\usepackage{systeme}
\usepackage{appendix}
\newenvironment{myproof}[1][\proofname]{\proof[#1]\mbox{}}{\endproof}

\usepackage[latin1]{inputenc}
\usepackage{setspace}

\newtheorem{theorem}{Theorem}[section]
\newtheorem{proposition}[theorem]{Proposition}

\newtheorem{cor}[theorem]{Corollary}
\newtheorem{observation}[theorem]{Remark}
\newtheorem{example}[theorem]{Example}
\begin{document}

\begin{center}
{\Large \bf{On Warped Product Gradient Yamabe Soliton}}

\bigskip
\bigskip
\bigskip
\textbf{\bf{Tokura, W. I.$^1$, Adriano, L. R.$^2$}, Pina R. S.$^3$.}\\
\textbf{\footnotesize \textit{ Instituto de Matemática e Estatística-UFG}}\\
\textbf{\footnotesize \textit{$^1$email:
williamisaotokura@hotmail.com}} \textbf{\footnotesize
\textit{$^2$email: levi@ufg.br}} \textbf{\footnotesize
\textit{$^3$email: romildo@ufg.br}}
\bigskip


{\Large}
\end{center}

    \begin{abstract}

In this paper, we provide a necessary and sufficient conditions for
the warped product $M=B\times_{f}F$ to be a gradient Yamabe soliton
when the base is conformal to an $n$-dimensional pseudo-Euclidean
space, which are invariant under the action of an
$(n-1)$-dimensional translation group, and the fiber $F$ is
scalar-constant. As application, we obtain solutions in steady case
with fiber scalar-flat. Besides, on the warped product we consider
the potential function as separable variables and obtain some
characterization of the base and the fiber.
\end{abstract}


\section{Introduction and main statements}

A \textit{Yamabe soliton} is a pseudo-Riemannian manifold $(M,g)$
admitting a vector field \linebreak$X\in \mathfrak{X}(M)$ such that

\begin{equation}\label{eq:01}(S_{g}-\rho)g=\frac{1}{2}\mathfrak{L}_{X}g,
\end{equation}
where $S_{g}$ denotes the scalar curvature of $M$, $\rho$ is a real
number and $\mathfrak{L}_{X}g$ denotes the Lie derivative of the
metric $g$ with respect to $X$. We say that $(M^{n},g)$ is
shrinking, steady or expanding, if $\rho>0$, $\rho=0$ , $\rho<0$,
respectively. When $X=\nabla h$ for some smooth function $h\in
\mathbf{C^{\infty}}(M)$, we say that $(M^{n},g,\nabla h)$ is an
\textit{gradient Yamabe soliton} with potential function $h$. In
this case the equation \eqref{eq:01} turns out

\begin{equation}\label{eq8}(S_{g}-\rho)g=Hess(h),
\end{equation}
where $Hess(\phi)$ denote the hessian of $\phi$. When $\phi$ is constant, we
call it a \textit{trivial Yamabe soliton}.

Yamabe solitons are self-similar solutions for the Yamabe flow
$$
\frac{\partial}{\partial t}g(t)=-R_{g(t)}g(t),
$$
and are important to
understand the geometric flow since they can appear as singularity
models. It has been known that every compact gradient Yamabe soliton
is of constant scalar curvature, hence, trivial since $f$ is
harmonic, see \cite{Cho}, \cite{Dask}, \cite{Hsu}. For the
non-compact case many interesting results is obtained in \cite{Cao}, \cite{Cat}, \cite{Ma}, \cite{Ma1}, \cite{Wu}.

As pointed in \cite{Cal} and \cite{Lopes}, it is important to
emphasize here that although the Yamabe flow are wellposed in the
Riemannian setting, they do not necessarily exist in the
semi-Riemannian case, where even the existence of short-time
solutions is not guaranteed in general due to the lack of
parabolicity. However, the existence of self-similar solutions of
the flow is equivalent to the existence of Yamabe solitons as in
\eqref{eq:01}. Semi-Riemannian Yamabe solitons have been intensively
studied, showing many differences with respect to the Riemannian
case, see for instance \cite{Bat} and \cite{Cal}.

Brozos-Vázquez et al. in \cite{Bro}, obtain a local characterization
of pseudo Riemannian ma\-ni\-fold endowed with gradient Yamabe soliton
metric, its results establish that if a pseudo Riemannian gradient
Yamabe soliton $(M,g)$, with potential function $h$ and such that
\linebreak$|\nabla h|\neq0$ is locally isometric to warped product of
unidimensional base and scalar-constant fiber. In the Riemannian
context a global structure result was given in \cite{Cao}.

In \cite{Dask} Daskalopoulos and Sesum investigated gradient Yamabe
soliton and proved that all complete locally conformally flat
gradient Yamabe solitons with positive sectional curvature are
rotationally symmetric. Proceeding in the same locally conformally
flat context, Neto and Tenenblat in \cite{Bene} consider the study of
pseudo Riemannian manifold
$(\mathbb{R}^{n},\frac{1}{\varphi^{}2}g_{0})$, where $g_{0}$ is the
canonical pseudo metric, and obtain a necessary and sufficient
condition to this manifold be a gradient Yamabe soliton. In the
search for invariant solutions they consider the invariant action of
an $(n-1)$-dimensional translation group and exhibit a complete
solution in steady case.

Recently, Pina and De Sousa in \cite{Pina}, consider the study of
gradient Ricci solitons on warped product structure
$M=B^{n}\times_{f}F^{m}$, where the base is conformal to an $n$-dimensional
pseudo-Euclidean space, invariant under the action of an
$(n-1)$-dimensional translation group, the fiber chosen to be an
Einstein manifold and potential function $h$ depending only on the
base and they give a necessary an sufficient condiction for $M$ to be a
gradient Ricci soliton.

As far as we know, there are no results for gradient Yamabe solitons
related to its potential function in the warped products of two
Riemannian manifolds of arbitrary dimensions. Thus, in this paper we
consider the study of gradient Yamabe solitons with warped product
structure, where we choose the fiber with dimension greater than
$1$. Initially we provide a sufficient condiction for the potential
function on warped product depends only on the base.

\begin{proposition}\label{eq3}Let $M=B\times_{f}F$ be a warped product
manifold with metric $\tilde{g}$. If the metric
$\tilde{g}=g_{B}\oplus f^{2}g_{F}$ is a gradient Yamabe soliton with
potential function $h:M\rightarrow \mathbb{R}$ and there exist a
pair of orthogonal vectors $(X_i,X_{j})$ of the base $B$, such that
$Hess_{g_{B}}(f)(X_{i},X_{j})\neq 0$, then the potential function
$h$ depends only on the base.
\end{proposition}

\begin{observation}
The previous proposition extend for Yamabe solitons the result obtained in \cite{Kim} where the authors studied warped product gradient Ricci solitons with one-dimensional base.
\end{observation}

Motivated by natural extension of Ricci solitons given by Rigoli,
Pigola and Setti in \cite{Pigola} the authors Barbosa and Ribeiro in
\cite{Barbosa}, define the concept of \textit{Almost Yamabe soliton}
allowing the constant $\rho$ in definition of Yamabe soliton
\eqref{eq:01} to be a differentiable function on $M$. The following
example was obtained by Barbosa and Ribeiro in \cite{Barbosa}, where
the manifold is endowed with a warped product metric.

\begin{example}Let $M^{n+1}=\mathbb{R}\times_{\cosh t}\mathbb{S}^{n}$
with metric $g=dt^{2}+\cosh^{2} t g_{0}$, where $g_{0}$ is the
canonical metric of $\mathbb{S}^{n}$. Taking $(M^{n+1},g,\nabla h,\rho)$, where $h(t,x)=\sinh t$ and \newline$\rho(t,x)=\sinh t+n$. A
straightforward computation gives that
$M^{n+1}=\mathbb{R}\times_{\cosh t}\mathbb{S}^{n}$ is a noncompact almost gradient Yamabe soliton.
\end{example}

In what follows, inspired by Proposition \ref{eq3}, we consider a
warped product gradient Yamabe soliton $M=B\times_{f}F$ with
potential function $h$ splitting of the form
\begin{equation}
\label{h12}
h(x,y)=h_{1}(x)+h_{2}(y),\ \mbox{where}\ h_{1}\in\mathcal{C}^{\infty}(B)\ 
\mbox{and}\ h_{2}\in\mathcal{C}^{\infty}(F),
\end{equation}
and get the following characterization theorem.

\begin{theorem}
\label{theo1.1}Let $M=B\times_{f}F$ be a warped product manifold with metric $\tilde{g}=g_{B}\oplus
f^{2}g_{F}$, and gradient Yamabe soliton structure with potential
function $h:B\times F\rightarrow\mathbb{R}$ given by \eqref{h12}, then one of the following
cases occurs

\begin{description}
  \item[(a)]$M$ is the Riemannian product between a \textbf{trivial gradiente Yamabe soliton}  and a \textbf{gradient Yamabe soliton}.
  \item[(b)]$M$ is the Riemannian product between two \textbf{gradient Yamabe solitons}.
  \item[(c)]$M$ is the warped product between a \textbf{Almost gradient Yamabe solitons} and a \textbf{trivial gradiente Yamabe soliton}.
\end{description}
\end{theorem}

This characterization theorem shows us that if we take the potential
function depending only on the base then the fiber $F$ is of
constant scalar curvature. In what follows we will take a warped
product gradient Yamabe soliton with potential function of the form
$h(x,y)=h_{1}(x)+constant$, the base conformal to an $n$-dimensional
pseudo-Euclidean space, and the fiber chosen to be an
scalar-constant space. More precisely, let $(\mathbb{R}^{n},g)$ be
the pseudo-Euclidean space, $n\geq3$ with coordinates
$x=(x_{1},\dots,x_{n})$ and $g_{ij}=\delta_{ij}\epsilon_{i}$ and let
$M=(\mathbb{R}^{n},\bar{g})\times_{f}F^{m}$ be a warped product
where $\bar{g}=\frac{1}{\varphi^{2}}g$, $F$ a semi-Riemannian
scalar-constant manifold with curvature $\lambda_{F}$, $m\geq1$,
$f$,$\varphi$, $h:\mathbb{R}^{n}\rightarrow\mathbb{R}$, smooth
functions, and $f$ is a positive function. Then we obtain necessary
and sufficient conditions  for the warped product metric
$g_{B}\oplus f^{2}g_{F}$ to be a gradient Yamabe soliton.

\begin{theorem}\label{eq:02}Let $(\mathbb{R}^{n},g)$ be a
pseudo-Euclidean space, $n\geq3$ with coordinates \newline$x=(x_{1},\dots,x_{n})$ and $g_{ij}=\delta_{ij}\epsilon_{i}$, and
let $M=(\mathbb{R}^{n},\bar{g})\times_{f}F^{m}$ be a warped product
where $\bar{g}=\frac{1}{\varphi^{2}}g$, $F$ a semi-Riemannian
scalar-constant manifold with curvature $\lambda_{F}$, $m\geq1$,
$f$,$\varphi$, $h:\mathbb{R}^{n}\rightarrow\mathbb{R}$, smooth
functions, and $f$ is a positive function. Then the warped product
metric $\tilde{g}$ is a gradient Yamabe soliton with potential
function $h$ if, and only if, the functions $f, \varphi, h$
satisfy

\begin{equation}\label{eq:19}h_{,x_{i}x_{j}}+\frac{\varphi_{,x_{j}}}{\varphi}h_{,x_{i}}+\frac{\varphi_{,x_{i}}}{\varphi}h_{,x_{j}}=0\hspace{1cm}
i\neq j,
\end{equation}

\begin{equation}\label{eq:20}
\begin{split}
\Big{[}(n-1)\left(2\varphi\sum_{k}\varepsilon_{k}\varphi_{,x_{k}x_{k}}-n\sum_{k}\varepsilon_{k}\varphi_{,x_{k}}^{2}\right)+\frac{\lambda_{F}}{f^{2}}-\frac{2m}{f}\left(\varphi^{2}\sum_{k}\varepsilon_{k}f_{,x_{k}x_{k}}-(n-2)\varphi\sum_{k}\varepsilon_{k}\varphi_{,x_{k}}f_{,x_{k}}\right)+\\
-\frac{m(m-1)}{f^2}\varphi^{2}\sum_{k}\varepsilon_{k}f_{,x_{k}}^{2}-\rho\Big{]}\frac{\varepsilon_{i}}{\varphi^{2}}=h_{,x_{i}x_{i}}+2\frac{\varphi_{,x_{i}}}{\varphi}h_{,x_{i}}-\varepsilon_{i}\sum_{k}\varepsilon_{k}\frac{\varphi_{,x_{k}}}{\varphi}h_{,x_{k}}\hspace{1cm}i=j,
\end{split}
\end{equation}

\begin{equation}\label{eq:21}
\begin{split}
(n-1)\left(2\varphi\sum_{k}\varepsilon_{k}\varphi_{,x_{k}x_{k}}-n\sum_{k}\varepsilon_{k}\varphi_{,x_{k}}^{2}\right)+\frac{\lambda_{F}}{f^{2}}-\frac{2m}{f}\left(\varphi^{2}\sum_{k}\varepsilon_{k}f_{,x_{k}x_{k}}-(n-2)\varphi\sum_{k}\varepsilon_{k}\varphi_{,x_{k}}f_{,x_{k}}\right)+\\
-\frac{m(m-1)}{f^2}\varphi^{2}\sum_{k}\varepsilon_{k}f_{,x_{k}}^{2}-\rho=\frac{\varphi^{2}}{f}\sum_{k}\varepsilon_{k}f_{,x_{k}}h_{,x_{k}}.
\end{split}
\end{equation}
\end{theorem}

In order to obtain solutions for equations in Theorem \ref{eq:02},
we consider $f$, $\varphi$ and $h$ invariant under the action of an
$(n-1)$-dimensional translation group, and
$\xi=\sum_{i=1}^{n}\alpha_{i}x_{i}, \alpha_{i}\in\mathbb{R}$, be a
basic invariant for the $(n-1)$-dimensional translation group, then
we obtain

\begin{theorem}\label{eq4}Let $(\mathbb{R}^{n},g)$ be a pseudo-Euclidean space, $n\geq3$ with coordinates \newline$x=(x_{1},\dots,x_{n})$, $g_{ij}=\delta_{ij}\epsilon_{i}$ and
let $M=(\mathbb{R}^{n},\bar{g})\times_{f}F^{m}$ be a warped product
where $\bar{g}=\frac{1}{\varphi^{2}}g$, $F$ a semi-Riemannian
scalar-constant manifold with curvature $\lambda_{F}$, $m\geq1$,
$f$,$\varphi$, $h:\mathbb{R}^{n}\rightarrow\mathbb{R}$, smooth
functions and $f>0$. Consider the functions $f(\xi)$,
$\varphi(\xi)$ and $h(\xi)$, where \newline$\xi=\sum_{k=1}^{n}\alpha_{k}x_{k},\alpha_{k}\in\mathbb{R}$ and
$\sum_{k=1}^{n}\varepsilon_{k}\alpha_{k}^{2}=\varepsilon_{k_{0}}$ or
$\sum_{k=1}^{n}\varepsilon_{k}\alpha_{k}^{2}=0$. Then the warped
product metric $\tilde{g}$ is a gradient Yamabe soliton with
potential function $h$ if, and only if, $f$, $h$ and $\varphi$,
satisfy

\begin{equation}\label{eq:09}
               \vspace{0,6cm}h''+2\frac{\varphi'h'}{\varphi}=0
\end{equation}

\begin{equation}\label{eq1}
\begin{split}
\varepsilon_{k_{0}}[(n-1)(2\varphi\varphi''-n(\varphi')^{2})-2\frac{m}{f}(\varphi^{2}f''-(n-2)\varphi\varphi'f')-&\frac{m(m-1)}{f^2}\varphi^{2}(f')^{2}+
\varphi'h'\varphi]\\&=\rho-\frac{\lambda_{F}}{f^{2}}
\end{split}
\end{equation}

\begin{equation}\label{eq2}
\begin{split}
\varepsilon_{k_{0}}[(n-1)(2\varphi\varphi''-n(\varphi')^{2})-2\frac{m}{f}(\varphi^{2}f''-(n-2)\varphi\varphi'f')-&\frac{m(m-1)}{f^2}\varphi^{2}(f')^{2}-
\frac{\varphi^{2}}{f}f'h']\\&=\rho-\frac{\lambda_{F}}{f^{2}}
\end{split}
\end{equation}
when
$\sum_{k=1}^{n}\varepsilon_{k}\alpha_{k}^{2}=\varepsilon_{k_{0}}$.

And

\begin{equation}\label{eq:10}
               \vspace{0,6cm}h''+2\frac{\varphi'h'}{\varphi}=0
\end{equation}
\begin{equation}\label{eq10}
               \rho-\frac{\lambda_{F}}{f^{2}}=0
\end{equation}
when $\sum_{k=1}^{n}\varepsilon_{k}\alpha_{k}^{2}=0.$
\end{theorem}

It is interesting to know how geometry of the fiber manifold $F$
affects the geometry of the warped product
$M=(\mathbb{R}^{n},\bar{g})\times_{f}F^{m}$. He in \cite{He0} has
shown that any complete steady gradient Yamabe soliton on
$\mathbb{R}\times_{f}F$ is necessarily isometric to the Riemannian
product with constant $f$ and $F$ being of zero scalar curvature.
Moreover, he showed that there is no complete steady gradient Yamabe
soliton on $\mathbb{R}\times_{f}F^{n}$ with $n\geq2$ and $F$
positive scalar constant manifold.

As consequence of Theorem \ref{eq4} in the context of lightlike
vector invariance and scalar-constant fiber, we prove that if $F$
has positive scalar constant curvature then there is no shirinking
or steady gradient Yamabe soliton
$M=(\mathbb{R}^{n},\bar{g})\times_{f}F^{m}$ and when $F$ has
negative constant scalar curvature, there is no expanding or steady
gradient Yamabe soliton \newline$M=(\mathbb{R}^{n},\bar{g})\times_{f}F^{m}$,
this is translated into the following corollary.

\begin{cor}\label{cor1.7}In the context of Theorem \ref{eq4}, if $X=\sum_{k}\alpha_{k}\frac{\partial}{\partial_{k}}$ is a lightlike vector, assume that $\lambda_{F}>0$, then there is no expanding or
steady gradient Yamabe soliton with warped metric $\tilde{g}$ and
potential function $h$. Similarly, if we assume that $\lambda_{F}<0$, then there is no shrinking or
steady gradient Yamabe soliton with warped metric $\tilde{g}$ and
potential function $h$.
\end{cor}

Now, by equations \eqref{eq:09} and \eqref{eq:10} in Theorem \ref{eq4} we easily see that a necessary condition for
$M=(\mathbb{R}^{n},\overline{g})\times_{f}F^{m}$ be a gradient
Yamabe soliton with invariant solution $f(\xi)$,
$\varphi(\xi)$ and $h(\xi)$, where
$\xi=\sum_{k=1}^{n}\alpha_{k}x_{k},\alpha_{k}\in\mathbb{R}$ is that
$h$ is a monotone function. That is,

\begin{equation}h'(\xi)=\frac{\alpha}{\varphi^{2}(\xi)},\nonumber
\end{equation}
 for some $\alpha\in\mathbb{R}$.

We provide solutions for ODE in Theorem \ref{eq4} in two cases:
$h'=0$ and $h'\neq0$, with metric $\tilde{g}=g_{B}\oplus f^{2}g_{F}$
be a steady gradient Yamabe soliton, i.e. $\rho=0$, and $F$ is a
scalar-flat pseudo-Riemannian manifold.

\begin{theorem}\label{eq5}In the context of Theorem \ref{eq4}, if $\sum_{k=1}^{n}\varepsilon_{k}\alpha_{k}^{2}=\varepsilon_{k_{0}}\neq 0$ and $F$ scalar-flat fiber, then the warped product metric $\tilde{g}$ is a steady gradient
Yamabe soliton with potential function $h$ and $h'\neq0$ if, and
only if, $f$, $h$ and $\varphi$, satisfies

\begin{equation}\label{eqa}f(\xi)=\frac{e^{c}}{\varphi(\xi)},
\end{equation}

\begin{equation}\label{eqb}h(\xi)=\alpha\int\frac{1}{\varphi^{2}(\xi)}d\xi,
\end{equation}

\begin{equation}\label{eqc}(n+m-1)(n+m+2)\int\frac{\varphi d\varphi}{\alpha-\frac{2\beta(n+m-1)(n+m+2)}{n+m-2}\varphi^{\frac{n+m}{2}+1}}=\xi+\nu,\hspace{0,2cm}c\in\mathbb{R}
\end{equation}
where $c$, $\nu$, $\alpha$, $\beta\in \mathbb{R}$ with
$\alpha\neq0$.
\end{theorem}

\begin{theorem}\label{eq7} In the context of Theorem \ref{eq4}, if $\sum_{k=1}^{n}\varepsilon_{k}\alpha_{k}^{2}=\varepsilon_{k_{0}}\neq 0$ and $F$ scalar-flat fiber, then, given a smooth function $\varphi>0$, the warped product metric
$\tilde{g}$ is a steady gradient Yamabe soliton with potential
function $h$ and $h'=0$ if, and only if, $f$ and $h$, satisfies

\begin{equation}\label{eq23}h(\xi)=\text{constant},
\end{equation}

\begin{equation}\label{eq22}f=\varphi^{\frac{n-2}{m+1}}e^{\int z_{p}d\xi}\left(\int e^{-(m+1)\int z_{p}d\xi}d\xi+\frac{2}{m+1}C\right)^{\frac{2}{m+1}}
\end{equation}
for an appropriate function $z_{p}$.

\end{theorem}

In the null case $\sum_{k=1}^{n}\varepsilon_{k}\alpha_{k}^{2}=0$ we
obtain

\begin{theorem}\label{eq6} In the context of Theorem \ref{eq4}, if $\sum_{k=1}^{n}\varepsilon_{k}\alpha_{k}^{2}=0$ and $F$ scalar-flat fiber, then, given smooth
functions $\varphi(\xi)$ and $f(\xi)$,the warped product metric
$\tilde{g}$ is a steady gradient Yamabe soliton with potential
function $h$ if, and only if

\begin{equation}h(\xi)=\alpha\int\frac{1}{\varphi^{2}(\xi)}d\xi.\nonumber
\end{equation}

\end{theorem}

\begin{observation}As pointed in \cite{Chenn}, see Theorem 3.6, a necessary and
sufficient condition to warped product $B\times_{f}F$ be a
conformally flat is that the function $f$ defines a global conformal
deformation such that $(B,\frac{1}{f^{2}}g_{B})$ is a space of
constant curvature $c$ and $F$ has constant curvature $-c$. With
this observation, we see that the solutions of theorems \ref{eq5},
\ref{eq7} and \ref{eq6} defines a non locally conformally flat
metric if the warping function $f$ is not constant.
\end{observation}

\begin{observation}As we can see in the proof of Theorem \ref{eq4}, if $\rho$ is a function defined only on the
base, then we can easily extend Theorem \ref{eq4} into context of
almost gradient Yamabe solitons. In the particular case of lightlike
vectors there are infinitely many solutions, that is, given
$\varphi$ and $f$

\begin{equation}\rho(\xi)=\frac{\lambda_{F}}{f(\xi)^{2}}\nonumber
\end{equation}

\begin{equation}h(\xi)=\alpha\int\frac{1}{\varphi^{2}(\xi)}d\xi\nonumber
\end{equation}
provide a family of almost gradient Yamabe soliton with warped
product structure.

\end{observation}

Before proving our main results, we present some examples
illustrating the above theorems.

\begin{example}In Theorem \ref{eq5}, consider $\beta=0$, then we
have

\begin{equation}f(\xi)=\frac{e^{c}\sqrt{(n-1)(n+m+2)}}{\sqrt{2\alpha(\xi+\nu)}},\nonumber
\end{equation}

\begin{equation}h(\xi)=\frac{\alpha}{2}(n-1)(n+m+2)\ln|\xi+\nu|\nonumber
\end{equation}

\begin{equation}\varphi(\xi)=\sqrt{\frac{2\alpha(\xi+\nu)}{(n-1)(n+m+2)}}\nonumber
\end{equation}
where $\alpha(\xi+\nu)>0$ and $c\in\mathbb{R}$. Thus, the metric
$\tilde{g}=\frac{1}{\varphi^{2}}g_{0}\oplus f^{2}g_{F}$ is a steady
gradient Yamabe soliton defined in the semi-space of Euclidean space
$\mathbb{R}^{n}$ with potential function $h$.

\end{example}

\begin{example}In Theorem \ref{eq7} consider the warped product $M=(R^{2},\bar{g})\times_{f}F^{3}$. Given the function $\varphi(\xi)=e^{\frac{3\xi^{2}}{4}}$
we have that

\begin{equation}\bar{g}=e^{-\frac{3\xi^{2}}{2}}g_0, \hspace{0,3cm} h(\xi)=\text{constant},\hspace{0,3cm}f(\xi)=e^{\frac{\xi}{2}}\nonumber
\end{equation}
defines a steady gradient Yamabe soliton in warped metric.

\end{example}

\begin{example}In Theorem \ref{eq6} consider the Lorentzian space $(\mathbb{R}^{n},g)$ with coordinates $(x_{1},\dots,x_{n})$
and signature $\epsilon_{1}=-1$, $\epsilon_{k}=1$ for all $k\geq2$,
and $F^{m}$ a complete scalar flat manifold. Let $\xi=x_{1}+x_{2}$
and choose $\varphi(\xi)=\frac{1}{1+\xi^{2}}$. Then, for
$\alpha\neq0$

\begin{equation}\bar{g}=(1+\xi^{2})^{2}g, \hspace{0,3cm} h(\xi)=\alpha\left(\xi+\frac{2}{3}\xi^{3}+\frac{1}{5}\xi^{5}\right),\hspace{0,3cm}f\in\mathcal{C}^{\infty}\nonumber
\end{equation}
defines a steady gradient Yamabe soliton
$(\mathbb{R}^{n},\bar{g})\times_{f}(F^{m},g_{flat})$ with potential
function $h$ and warping function $f$. Observe that, since the
conformal function $\varphi$ is bounded, we have that $\bar{g}$ is
complete, and consequently
$\tilde{g}=\frac{1}{\varphi^{2}}g_{0}\oplus f^{2}g_{F}$ is complete.

\end{example}


\begin{section}{Proofs of the Main Results}

\begin{myproof}[\textbf{Proof of Proposition 1.1}]Let $M=B\times_{f}F$ be a gradient Yamabe soliton with potential function $h:M\rightarrow\mathbb{R}$,
then by equation \eqref{eq8}, we obtain

\begin{equation}\label{eq9}(S_{\tilde{g}}-\rho)\tilde{g}=Hess(h)
\end{equation}

Now, It is well known that on warped metric $\tilde{g}$ the scalar
curvature of base $B$, fiber $F$ and $M$ is related by(see chapter 7
of \cite{oneil})

\begin{equation}S_{\tilde{g}}=S_{g_{B}}+\frac{S_{g_{F}}}{f^{2}}-\frac{2d}{f}\Delta^{B}f-d(d-1)\frac{\langle grad_{B}f,grad_{B}f\rangle}{f^{2}}\label{eq:04}
\end{equation}
where $\Delta^{B}$ denote the laplacian defined on $B$. Then
considering $X_{1},X_{2},\dots,X_{n}\in \mathcal{L}(B)$ and
$Y_{1},Y_{2},\dots,Y_{m}\in \mathcal{L}(F)$, where $\mathcal{L}(B)$
and $\mathcal{L}(F)$ are respectively the space of lifts of vector
fiels on $B$ and $F$ to $B\times F$, we obtain substituting equation
\eqref{eq:04} in \eqref{eq9} that

\begin{equation*}
\begin{cases}
               \vspace{0,6cm}
               \Big{(}S_{g_{B}}+\frac{S_{g_{F}}}{f^{2}}-\frac{2d}{f}\Delta^{B}f-d(d-1)\frac{\langle
grad_{B}f,grad_{B}f\rangle}{f^{2}}-\rho\Big{)}g_{B}(X_{i},X_{j})=Hess_{\tilde{g}}h(X_{i},X_{j})&(i)\\
               \vspace{0,6cm}
               \Big{(}S_{g_{B}}+\frac{S_{g_{F}}}{f^{2}}-\frac{2d}{f}\Delta^{B}f-d(d-1)\frac{\langle
grad_{B}f,grad_{B}f\rangle}{f^{2}}-\rho\Big{)}\tilde{g}(X_{i},Y_{j})=Hess_{\tilde{g}}h(X_{i},Y_{j})&(ii)  \\

               \Big{(}S_{g_{B}}+\frac{S_{g_{F}}}{f^{2}}-\frac{2d}{f}\Delta^{B}f-d(d-1)\frac{\langle
grad_{B}f,grad_{B}f\rangle}{f^{2}}-\rho\Big{)}f^{2}g_{F}(Y_{i},Y_{j})=Hess_{\tilde{g}}h(Y_{i},Y_{j}).&(iii)
\end{cases}
\end{equation*}
Thus, using the fact $\tilde{g}(X_{i},Y_{j})=0$, we obtain by expression $(ii)$ that

\begin{equation}Hess(h)(X_{i},Y_{j})=0.\nonumber
\end{equation}

Now, using Lemma 2.1 of \cite{He}, we obtain

\begin{equation}\label{eqe}h(x,y)=z(x)+f(x)v(y)
\end{equation}
where $z:B\rightarrow\mathbb{R}$ and $v:F\rightarrow\mathbb{R}$.
Then since

$$Hess_{\tilde{g}}h(X_{i},X_{j})=Hess_{g_{B}}h(X_{i},X_{j}),$$
we have by expression $(i)$ that

\begin{equation}\label{eqd}\Big{(}S_{g_{B}}+\frac{S_{g_{F}}}{f^{2}}-\frac{2d}{f}\Delta^{B}f-d(d-1)\frac{\langle
grad_{B}f,grad_{B}f\rangle}{f^{2}}-\rho\Big{)}g_{B}(X_{i},X_{j})=Hess_{g_{B}}h(X_{i},X_{j}).
\end{equation}
We have from \eqref{eqe} that the right size of \eqref{eqd} is given by

\begin{equation}\label{eqf}Hess_{g_{B}}z+vHess_{g_{B}}f.
\end{equation}
By hypothesis, there are two ortogonal vectors fields $X_{i}$, $X_{j}$ such that $Hess_{g_{B}}f(X_{i},X_{j})\neq0$. Then combining this fact with equations \eqref{eqd} and \eqref{eqf}, we obtain

\begin{equation}v=-\frac{Hess_{g_{B}}z(X_{i},X_{j})}{Hess_{g_{B}}f(X_{i},X_{j})}.
\end{equation}
This show that $v(y)$ is constant, and then by expression \eqref{eqe} we have that $h$ depends only on the base.

\end{myproof}

\begin{myproof}[\textbf{Proof of Theorem \ref{theo1.1}}]Let $M=B\times_{f}F$ be a warped product with gradient Yamabe soliton structure and potential function
 $h(x,y)=h_{1}(x)+h_{2}(y)$. In the same way as in the proof of Proposition \ref{eq3},
 for $X_{1},X_{2},\dots,X_{n}\in \mathcal{L}(B)$ and
$Y_{1},Y_{2},\dots,Y_{m}\in \mathcal{L}(F)$ we obtain

\begin{equation}Hess(h)(X_{i},Y_{j})=0.\nonumber
\end{equation}

As we know, the connection of warped product is particularly simple, that is, for \newline$X\in \mathcal{L}(B)$ and $Y\in \mathcal{L}(F)$, we have

$$\nabla_{X}Y=\nabla_{Y}X=\frac{X(f)}{f}Y.$$
Thus,

\begin{equation}Hess(h)(X_{i},Y_{j})=X_{i}(Y_{j}(h))-(\nabla_{X_{i}}Y_{j})h=X_{i}(Y_{j}(h))-\frac{X_{i}(f)}{f}Y_{j}=0.\nonumber
\end{equation}

Establishing the notation $h_{,x_{i}}=X_{i}(h)$,
$h_{,x_{i}x_{j}}=X_{j}(X_{i}(h))$, we have that

\begin{equation}h_{,y_{j}x_{i}}-\frac{f_{,x_{i}}}{f}h_{,y_{j}}=0-\frac{f_{,x_{i}}}{f}(h_{2})_{,y_{j}}=0\hspace{0,4cm}\forall i,j.\nonumber
\end{equation}
Then, $f$ is constant or $h(x,y)=h_{1}(x)+constant$. We separate the
proof in three cases:

$Case (I):$ ($f$ is constant and $h(x,y)=h_{1}(x)+constant$). In
this case, $M=B\times_{f}F$ is a Riemannian product and we have

\begin{equation*}
\begin{cases}
               \vspace{0,6cm}
               \Big{(}S_{g_{B}}+\frac{S_{g_{F}}}{f^{2}}-\rho\Big{)}g_{B}(X_{i},X_{j})=Hess_{g_{B}}h_{1}(X_{i},X_{j})&(i)\\
               \vspace{0,6cm}
               \Big{(}S_{g_{B}}+\frac{S_{g_{F}}}{f^{2}}-\rho\Big{)}\tilde{g}(X_{i},Y_{j})=Hess_{\tilde{g}}h(X_{i},Y_{j})=0&(ii)  \\

               \Big{(}S_{g_{B}}+\frac{S_{g_{F}}}{f^{2}}-\rho\Big{)}f^{2}g_{F}(Y_{i},Y_{j})=Hess_{g_{F}}h_{2}(Y_{i},Y_{j})+f\nabla f(h_{1})g_{F}(Y_{i},Y_{j})&(iii)
\end{cases}
\end{equation*}
where in equation $(iii)$ we use Proposition 35 of \cite{oneil} and
the Hessian definition to get

\begin{eqnarray}\label{eq13}
Hess_{\tilde{g}}h(Y_{i},Y_{j})& = &
Y_{i}(Y_{j}(h))-(\nabla_{Y_{i}}Y_{j})^{M}h\\
&=& Y_{i}(Y_{j}(h))-(\mathcal{H}(\nabla_{Y_{i}}Y_{j})+\mathcal{V}(\nabla_{Y_{i}}Y_{j}))(h)\nonumber\\
&=& Y_{i}(Y_{j}(h))+\frac{\langle Y_{i},
  Y_{j}\rangle}{f}grad_{\tilde{g}}f(h)-\nabla_{Y_{i}}^{F}Y_{j}(h)\nonumber\\
&=& Y_{i}(Y_{j}(h))+fg_{F}(Y_{i},
  Y_{j})grad_{\tilde{g}}f(h)-\nabla_{Y_{i}}^{F}Y_{j}(h)\nonumber\\
  &=&Hess_{g_{F}}h_{2}(Y_{i},Y_{j})+f\nabla
  f(h_{1})g_{F}(Y_{i},Y_{j}).\nonumber
\end{eqnarray}

Since  $S_{g_{F}}$ is constant on $B$, we have from $(i)$ that $B$
is a gradient Yamabe soliton of the form $(B,g_{B},\nabla
h_{1},-\frac{S_{g_{F}}}{f^{2}}+\rho)$. Furthermore, since
$h(x,y)=h_{1}(x)+cte$ we have by $(iii)$ that $F$ is a trivial
gradient Yamabe soliton of the form $(F,g_{F},\nabla
0,f^{2}\rho-f^{2}S_{g_{B}})$. This proves the item $(a)$.

\vspace{0,5cm}

$Case (II):$ ($f$ is constant and $h(x,y)=h_{1}(x)+h_{2}$, $h_{2}$
not necessarily constant). In this case, $M=B\times_{f}F$ is a
Riemannian product and we have

\begin{equation*}
\begin{cases}
               \vspace{0,6cm}
               \Big{(}S_{g_{B}}+\frac{S_{g_{F}}}{f^{2}}-\rho\Big{)}g_{B}(X_{i},X_{j})=Hess_{g_{B}}h_{1}(X_{i},X_{j})&(i)\\
               \vspace{0,6cm}
               \Big{(}S_{g_{B}}+\frac{S_{g_{F}}}{f^{2}}-\rho\Big{)}\tilde{g}(X_{i},Y_{j})=Hess_{\tilde{g}}h(X_{i},Y_{j})=0&(ii)  \\

               \Big{(}S_{g_{B}}+\frac{S_{g_{F}}}{f^{2}}-\rho\Big{)}f^{2}g_{F}(Y_{i},Y_{j})=Hess_{g_{F}}h_{2}(Y_{i},Y_{j})+f\nabla f(h_{1})g_{F}(Y_{i},Y_{j}).&(iii)
\end{cases}
\end{equation*}

Since $S_{g_{F}}$ is constant on $B$, we have that $B$ is a gradient
Yamabe soliton of the form $(B,g_{B},\nabla
h_{1},-\frac{S_{g_{F}}}{f^{2}}+\rho)$. Furthermore, by equation
$(iii)$ we have that $F$ is a gradient Yamabe soliton of the form
$(F,g_{F},\nabla h_{2},f^{2}\rho-f^{2}S_{g_{B}})$. This proves the
item $(b)$.

$Case (III):$ ($f$ is non constant and $h(x,y)=h_{1}(x)+constant$).
In this case we have:

\begin{equation*}
\begin{cases}
               \vspace{0,6cm}
               \Big{(}S_{g_{B}}+\frac{S_{g_{F}}}{f^{2}}-\frac{2d}{f}\Delta^{B}f-d(d-1)\frac{\langle
grad_{B}f,grad_{B}f\rangle}{f^{2}}-\rho\Big{)}g_{B}(X_{i},X_{j})=Hess_{g_{B}}h_{1}(X_{i},X_{j})&(i)\\
               \vspace{0,6cm}
               \Big{(}S_{g_{B}}+\frac{S_{g_{F}}}{f^{2}}-\frac{2d}{f}\Delta^{B}f-d(d-1)\frac{\langle
grad_{B}f,grad_{B}f\rangle}{f^{2}}-\rho\Big{)}\tilde{g}(X_{i},Y_{j})=Hess_{\tilde{g}}h(X_{i},Y_{j})=0&(ii)  \\

               \Big{(}S_{g_{B}}+\frac{S_{g_{F}}}{f^{2}}-\frac{2d}{f}\Delta^{B}f-d(d-1)\frac{\langle
grad_{B}f,grad_{B}f\rangle}{f^{2}}-\rho\Big{)}f^{2}g_{F}(Y_{i},Y_{j})=f\nabla
f(h)g_{F}(Y_{i},Y_{j})&(iii)
\end{cases}
\end{equation*}

Since $f>0$, by equation $(iii)$ we have that
\begin{equation} (S_{g_{F}}-\psi) g_{F}(Y_{i},Y_{j})=0
\end{equation}
where $\psi=-f^{2}S_{g_{B}}+2fd\Delta^{B}f+d(d-1)\langle
grad_{B}f,grad_{B}f\rangle+f^{2}\rho+f\nabla f(h)$.

Now, since $\psi$ depend only on $B$, we have that $\psi$ is
constant on $F$, then by equation $(iii)$, we have that $F$ is a
trivial gradient Yamabe soliton. Furthermore, by equation $(i)$ we
have that $(B,g_{B})$ is a gradient almost Yamabe soliton of the
form 
$$
(B,g_{B},\nabla
h_{1},-[\frac{S_{g_{F}}}{f^{2}}-\frac{2d}{f}\Delta^{B}f-d(d-1)\frac{\langle
grad_{B}f,grad_{B}f\rangle}{f^{2}}-\rho] ).
$$
This proves the item $(c)$.
\end{myproof}

\begin{myproof}[\textbf{Proof of Theorem \ref{eq:02}}]
Let $M$ be a warped product with gradient Yamabe soliton structure
and potential function $h$, that is,
\begin{equation}(S_{\tilde{g}}-\rho)\tilde{g}=Hess_{\tilde{g}}(h).\label{eq:200}
\end{equation}

By the same arguments used in proof of Proposition \ref{eq3}, for
$X_{1},X_{2},\dots,X_{n}\in \mathcal{L}(B)$ and
$Y_{1},Y_{2},\dots,Y_{m}\in \mathcal{L}(F)$ we obtain

\begin{equation*}
\begin{cases}
               \vspace{0,6cm}
               \Big{(}S_{\overline{g}}+\frac{S_{g_{F}}}{f^{2}}-\frac{2m}{f}\Delta_{\overline{g}}f-m(m-1)\frac{\langle
grad_{\overline{g}}f,grad_{\overline{g}}f\rangle}{f^{2}}-\rho\Big{)}\overline{g}(X_{i},X_{j})=Hess_{\tilde{g}}h(X_{i},X_{j})&(i)\\
               \vspace{0,6cm}
               \Big{(}S_{\overline{g}}+\frac{S_{g_{F}}}{f^{2}}-\frac{2m}{f}\Delta_{\overline{g}}f-m(m-1)\frac{\langle
grad_{\overline{g}}f,grad_{\overline{g}}f\rangle}{f^{2}}-\rho\Big{)}\tilde{g}(X_{i},Y_{j})=Hess_{\tilde{g}}h(X_{i},Y_{j})=0&(ii)  \\

               \Big{(}S_{\overline{g}}+\frac{S_{g_{F}}}{f^{2}}-\frac{2m}{f}\Delta_{\overline{g}}f-m(m-1)\frac{\langle
grad_{\overline{g}}f,grad_{\overline{g}}f\rangle}{f^{2}}-\rho\Big{)}f^{2}g_{F}(Y_{i},Y_{j})=Hess_{\tilde{g}}h(Y_{i},Y_{j}).&(iii)
\end{cases}
\end{equation*}

It is well known that for the conformal metric
$\bar{g}=\frac{1}{\varphi^{2}}g_{0}$, the Christofel symbol is given
by

\begin{equation}\bar{\Gamma}_{ij}^{k}=0,\ \bar{\Gamma}_{ij}^{i}=-\frac{\varphi_{,x_{j}}}{\varphi},\ \bar{\Gamma}_{ii}^{k}=\epsilon_{i}\epsilon_{k}\frac{\varphi_{,x_{k}}}{\varphi}\;\ \mbox{and}\;\ \bar{\Gamma}_{ii}^{i}=-\frac{\varphi_{,x_{i}}}{\varphi}.\nonumber
\end{equation}
Then, we obtain by Hessian definiton that
\begin{equation}\label{eq:22}
\begin{cases}
               \vspace{0,6cm}
               Hess_{\overline{g}}(h)_{ij}=h_{,x_{i}x_{j}}+\frac{\varphi_{,x_{i}}h_{,x_{j}}}{\varphi}+\frac{\varphi_{,x_{j}}h_{,x_{i}}}{\varphi}& i\neq j\\

               Hess_{\overline{g}}(h)_{ii}=h_{,x_{i}x_{i}}+2\frac{\varphi_{,x_{i}}h_{,x_{i}}}{\varphi}-\varepsilon_{i}\sum_{k=1}^{n}\varepsilon_{k}\frac{\varphi_{,x_{k}}}{\varphi}h_{,x_{k}}&i=j.
\end{cases}
\end{equation}
The Ricci curvature is given by

$$Ric_{\overline{g}}=\frac{1}{\varphi^{2}}\Big{\{}(n-2)\varphi
Hess_{g}(\varphi)+[\varphi\Delta_{g}\varphi-(n-1)|\nabla_{g}\varphi|^{2}]g\Big{\}}$$
and then we easily see that the scalar curvature on conformal metric
is given by

\begin{equation}\label{eq12}S_{\overline{g}}=(n-1)(2\varphi\Delta_{g}\varphi-n|\nabla_{g}\varphi|^{2})=(n-1)(2\varphi\sum_{k=1}^{n}\varepsilon_{k}\varphi_{,x_{k}x_{k}}-n\sum_{k=1}^{n}\varepsilon_{k}\varphi_{,x_{k}}^{2}).
\end{equation}
Since $h:\mathbb{R}^{n}\rightarrow\mathbb{R}$, we obtain

\begin{equation}\label{eq:24}
Hess_{\tilde{g}}h(X_{i},X_{j})=Hess_{\overline{g}}h(X_{i},X_{j}),
\hspace{0,4cm} \forall i,j.
\end{equation}

On the other hand

\begin{equation}\label{eq:23}
\begin{cases}
               \vspace{0,4cm}

               S_{F}g_{F}=\lambda_{F}g_{F}\\
\vspace{0,4cm}
               \tilde{g}(Y_{i},Y_{j})=f^{2}g_{F}(Y_{i},Y_{j})\\
\vspace{0,4cm}
               \Delta_{\overline{g}}f=\varphi^{2}\sum_{k}\varepsilon_{k}f_{,x_{k}x_{k}}-(n-2)\varphi\sum_{k}\varepsilon_{k}\varphi_{,x_{k}}f_{,x_{k}}\\

               \tilde{g}(grad_{\overline{g}}f,grad_{\overline{g}}f)=\overline{g}(grad_{\overline{g}}f,grad_{\overline{g}}f)=\varphi^{2}\sum_{k}\varepsilon_{k}f_{,x_{k}}^{2}.

\end{cases}
\end{equation}

Now, substituting the second equation of $\eqref{eq:22}$, the
equations of \eqref{eq:23} and equation \eqref{eq12} in $(i)$, we
have
\begin{equation*}
\begin{split}
\Big{[}(n-1)(2\varphi\sum_{k}\varepsilon_{k}\varphi_{,x_{k}x_{k}}-n\sum_{k}\varepsilon_{k}\varphi_{,x_{k}}^{2})+\frac{\lambda_{F}}{f^{2}}-\frac{2m}{f}(\varphi^{2}\sum_{k}\varepsilon_{k}f_{,x_{k}x_{k}}-(n-2)\varphi\sum_{k}\varepsilon_{k}\varphi_{,x_{k}}f_{,x_{k}})+\\
-\frac{m(m-1)}{f^2}\varphi^{2}\sum_{k}\varepsilon_{k}f_{,x_{k}}^{2}-\rho\Big{]}\frac{\varepsilon_{i}}{\varphi^{2}}=h_{,x_{i}x_{i}}+2\frac{\varphi_{,x_{i}}}{\varphi}h_{,x_{i}}-\varepsilon_{i}\sum_{k}\varepsilon_{k}\frac{\varphi_{,x_{k}}}{\varphi}h_{,x_{k}}
\end{split}
\end{equation*}
which is the equation \eqref{eq:20}.

Analogously, substituting the first equation of \eqref{eq:22} and
equation \eqref{eq:23} in $(i)$, we obtain
\begin{equation}h_{,x_{i}x_{j}}+\frac{\varphi_{,x_{j}}}{\varphi}h_{,x_{i}}+\frac{\varphi_{,x_{i}}}{\varphi}h_{,x_{j}}=0
\end{equation}
which is the equation \eqref{eq:19}.

in the similar way that equation \eqref{eq13}, we have that

\begin{eqnarray}\label{eq:25}
Hess_{\tilde{g}}h(Y_{i},Y_{j})& = &
Y_{i}(Y_{j}(h))-(\nabla_{Y_{i}}Y_{j})^{M}h\nonumber\\
 &=&
Hess_{gF}h(Y_{i},Y_{j})+fg_{F}(Y_{i},
  Y_{j})grad_{\tilde{g}}f(h)\nonumber\\
&=& fg_{F}(Y_{i},Y_{j})grad_{\tilde{g}}f(h)\nonumber\\
&=&f\varphi^{2}\sum_{k}\varepsilon_{k}f_{,x_{k}}h_{,x_{k}}g_{F}(Y_{i},Y_{j}).
\end{eqnarray}
Then, substituting equation \eqref{eq:25}, \eqref{eq:23} and
equation \eqref{eq12} in $(iii)$, we obtain equation \eqref{eq:21}.
A direct calculation shows us the converse implication. This
concludes the proof of Theorem \ref{eq:02}.
\end{myproof}

\begin{myproof}[\textbf{Proof of Theorem \ref{eq4}}]Since we are assuming that $\varphi(\xi)$, $h(\xi)$
and $f(\xi)$ are functions of $\xi$, where
$\xi=\sum_{k}\alpha_{k}x_{k}$, $\alpha_{k}\in\mathbb{R}^{n}$ and
$\sum_{k=1}^{n}\varepsilon_{k}\alpha_{k}^{2}=\varepsilon_{k_{0}}$ or
$\sum_{k=1}^{n}\varepsilon_{k}\alpha_{k}^{2}=0$, then we have

\begin{equation*}\varphi_{,x_{i}}=\varphi'\alpha_{i};\hspace{0,1cm}
\varphi_{,x_{i}x_{j}}=\varphi''\alpha_{i}\alpha_{j};\hspace{0,1cm}f{,x_{i}}=f'\alpha_{i};\hspace{0,1cm}
f_{,x_{i}x_{j}}=f''\alpha_{i}\alpha_{j};\hspace{0,1cm}
h_{,x_{i}}=h''\alpha_{i};\hspace{0,1cm}
h_{,x_{i}x_{j}}=h''\alpha_{i}\alpha_{j}.
\end{equation*}
Substituting these expressions into \eqref{eq:19} of Theorem \ref{eq:02}, we obtain

\begin{equation}\label{eq:33}\Big{(}h''+2\frac{h'\varphi'}{\varphi}\Big{)}\alpha_{i}\alpha_{j}=0,
\hspace{0,4cm} \forall i\neq j.
\end{equation}
Similarly, considering equations \eqref{eq:20} and \eqref{eq:21} of Theorem \ref{eq:02},
we obtain

\begin{equation}\label{eq:30}
\begin{split}
\Big{[}(n-1)(2\varphi\varphi''\sum_{k}\varepsilon_{k}\alpha_{k}^{2}-n(\varphi')^{2}\sum_{k}\varepsilon_{k}\alpha_{k}^{2})+\frac{\lambda_{F}}{f^{2}}-\frac{2m}{f}(\varphi^{2}f''\sum_{k}\varepsilon_{k}\alpha_{k}^{2}-(n-2)\varphi\varphi'f'\sum_{k}\varepsilon_{k}\alpha_{k}^{2})+\\
-\frac{m(m-1)}{f^2}\varphi^{2}(f')^{2}\sum_{k}\varepsilon_{k}\alpha_{k}^{2}-\rho\Big{]}\frac{\varepsilon_{i}}{\varphi^{2}}=h''\alpha_{i}^{2}+2\alpha_{i}^{2}\frac{\varphi'}{\varphi}h'-\varepsilon_{i}h'\frac{\varphi'}{\varphi}\sum_{k}\varepsilon_{k}\alpha_{k}^{2}
\end{split}
\end{equation}
for $i\in\{1,2,\dots,n\}$, and

\begin{equation}\label{eq:31}
\begin{split}
(n-1)(2\varphi\varphi''\sum_{k}\varepsilon_{k}\alpha_{k}^{2}-n(\varphi')^{2}\sum_{k}\varepsilon_{k}\alpha_{k}^{2})+\frac{\lambda_{F}}{f^{2}}-\frac{2m}{f}(\varphi^{2}f''\sum_{k}\varepsilon_{k}\alpha_{k}^{2}-(n-2)\varphi\varphi'f'\sum_{k}\varepsilon_{k}\alpha_{k}^{2})+\\
-\frac{m(m-1)}{f^2}\varphi^{2}(f')^{2}\sum_{k}\varepsilon_{k}\alpha_{k}^{2}-\rho=\frac{\varphi^{2}}{f}f'h'\sum_{k}\varepsilon_{k}\alpha_{k}^{2}.
\end{split}
\end{equation}

If there exist $i,j$, $i\neq j$ such that $\alpha_{i}\alpha_{j}\neq
0$, then we get by equation \eqref{eq:33} that

\begin{equation}\label{eq:34}\left(h''+2\frac{h'\varphi'}{\varphi}\right)=0.
\end{equation}
It follows from \eqref{eq:34} that the equation \eqref{eq:30} is
summed to

\begin{equation}
\begin{split}
\Big{[}(n-1)(2\varphi\varphi''\sum_{k}\varepsilon_{k}\alpha_{k}^{2}-n(\varphi')^{2}\sum_{k}\varepsilon_{k}\alpha_{k}^{2})+\frac{\lambda_{F}}{f^{2}}-\frac{2m}{f}(\varphi^{2}f''\sum_{k}\varepsilon_{k}\alpha_{k}^{2}-(n-2)\varphi\varphi'f'\sum_{k}\varepsilon_{k}\alpha_{k}^{2})+\\
-\frac{m(m-1)}{f^2}\varphi^{2}(f')^{2}\sum_{k}\varepsilon_{k}\alpha_{k}^{2}-\rho\Big{]}\frac{\varepsilon_{i}}{\varphi^{2}}=-\varepsilon_{i}h'\frac{\varphi'}{\varphi}\sum_{k}\varepsilon_{k}\alpha_{k}^{2}.
\end{split}
\end{equation}
Thus, isolating the term
$\sum_{k}\varepsilon_{k}\alpha_{k}^{2}=\varepsilon_{k_{0}}$, we
obtain the equation \eqref{eq1}.

In the same way, isolating the term
$\sum_{k}\varepsilon_{k}\alpha_{k}^{2}=\varepsilon_{k_{0}}$ in
\eqref{eq:31} we obtain the equation \eqref{eq2}.

Thus, if
$\sum_{k}\varepsilon_{k}\alpha_{k}^{2}=\varepsilon_{k_{0}}$, then we
obtain equations \eqref{eq:09}, \eqref{eq1} and \eqref{eq2}. In the
case \newline$\sum_{k}\varepsilon_{k}\alpha_{k}^{2}=0$, we easily see that
the equation \eqref{eq:09} is summed to

\begin{equation}
\begin{cases}
               \vspace{0,6cm}h''+2\frac{\varphi'h'}{\varphi}=0 \\

               \rho-\frac{\lambda_{F}}{f^{2}}=0.\\
\end{cases}
\end{equation}

Now, we need to consider the case $\alpha_{k_{0}}=1$ and
$\alpha_{k}=0$ $\forall k\neq k_{0}$. In this case, equation
\eqref{eq:33} is trivially satisfied, and since equation
\eqref{eq:31} does not depend on the index $i$, we have that
equation \eqref{eq:31} is equivalent to equation $\eqref{eq2}$.

Finally, we need to show the validity of equation \eqref{eq:09} and
$\eqref{eq1}$. Observe that taking $i=k_{0}$, that is,
$\alpha_{k_{0}}=1$, in \eqref{eq:30}, we get

\begin{equation}
\begin{split}
\Big{[}(n-1)(2\varphi\varphi''\varepsilon_{k_{0}}-n(\varphi')^{2}\varepsilon_{k_{0}})+\frac{\lambda_{F}}{f^{2}}-\frac{2m}{f}(\varphi^{2}f''\varepsilon_{k_{0}}-(n-2)\varphi\varphi'f'\varepsilon_{k_{0}})+\\
-\frac{m(m-1)}{f^2}\varphi^{2}(f')^{2}\varepsilon_{k_{0}}-\rho\Big{]}\frac{\varepsilon_{k_{0}}}{\varphi^{2}}=h''+2\frac{\varphi'}{\varphi}h'-h'\frac{\varphi'}{\varphi}=h''+\frac{\varphi'}{\varphi}h'
\end{split}
\end{equation}
and for $i\neq k_{0}$, that is, $\alpha_{i}=0$, we have

\begin{equation}
\begin{split}
\Big{[}(n-1)(2\varphi\varphi''\varepsilon_{k_{0}}-n(\varphi')^{2}\varepsilon_{k_{0}})+\frac{\lambda_{F}}{f^{2}}-\frac{2m}{f}(\varphi^{2}f''\varepsilon_{k_{0}}-(n-2)\varphi\varphi'f'\varepsilon_{k_{0}})+\\
-\frac{m(m-1)}{f^2}\varphi^{2}(f')^{2}\varepsilon_{k_{0}}-\rho\Big{]}\frac{\varepsilon_{i}}{\varphi^{2}}=-\varepsilon_{i}\varepsilon_{k_{0}}\frac{\varphi'}{\varphi}h'.
\end{split}
\end{equation}

However, this equations are equivalent to equations \eqref{eq:09}
and $\eqref{eq1}$. This complete the proof of Theorem \ref{eq4}.
\end{myproof}

\begin{myproof}[\textbf{Proof of Corollary \ref{cor1.7}}]By Theorem \ref{eq4}, we have that $M$ is a gradient
Yamabe soliton with potential function $h$ if, and only if,

\begin{equation}
\begin{cases}
               \vspace{0,6cm}h''+2\frac{\varphi'h'}{\varphi}=0 \\

               \rho-\frac{\lambda_{F}}{f^{2}}=0.\\
\end{cases}
\end{equation}
Thus, we have that $\lambda_{F}$ and $\rho$ always
have the same signal. Therefore, there is no gradient Yamabe soliton
$M$ expanding/shrinking with fiber trivial gradient Yamabe soliton
shrinking/expanding.
\end{myproof}

\begin{myproof}[\textbf{Proof of Theorem \ref{eq5}}]Since $\lambda_{F}=\rho=0$ we have by equation \eqref{eq1} and \eqref{eq2}
of Theorem \ref{eq4} that

\begin{equation}\varphi'h'\varphi=-\frac{\varphi^{2}}{f}f'h'\nonumber
\end{equation}
and by condition $h'\neq0$, we obtain
\begin{equation}\label{eq14}\frac{\varphi'}{\varphi}=-\frac{f'}{f}.
\end{equation}
Integrating this equation we have

$$f(\xi)=\frac{e^{c}}{\varphi(\xi)}$$
for some $c\in\mathbb{R}$, which is equation \eqref{eqa} of Theorem \ref{eq5}.

Integrating the equation \eqref{eq:09}, we have that

\begin{equation}\label{eq19}h'(\xi)=\frac{\alpha}{\varphi^{2}(\xi)}
\end{equation}
for some $\alpha\neq0$, and
$$h(\xi)=\alpha\int\frac{1}{\varphi^{2}(\xi)}d\xi$$
which is equation \eqref{eqb} of Theorem \ref{eq5}.

Substituting equation \eqref{eq19} into \eqref{eq1}
we have

\begin{equation}\label{eq15}(n-1)(2\varphi\varphi''-n(\varphi')^{2})-2\frac{m}{f}(\varphi^{2}f''-(n-2)\varphi\varphi'f')-\frac{m(m-1)}{f^2}\varphi^{2}(f')^{2}+\alpha\frac{\varphi'}{\varphi}=0.
\end{equation}
Inserting equation \eqref{eq14} into \eqref{eq15} we obtain

\begin{equation}\label{eq16}\varphi\varphi''-\frac{(n+m)}{2}(\varphi')^{2}+\frac{\alpha}{2(n+m-1)}\frac{\varphi'}{\varphi}=0.
\end{equation}
Consider $\varphi(\xi)^{1-\frac{n+m}{2}}=\omega(\xi)$, then
\begin{equation}\omega'(\xi)=(1-\frac{n+m}{2})\varphi^{-\frac{n+m}{2}}\varphi',\hspace{0,5cm}\omega''(\xi)=(1-\frac{m+n}{2})\left(\varphi^{-\frac{n+m}{2}-1}(\varphi\varphi''-\frac{(n+m)}{2}(\varphi')^{2})\right)
\end{equation}
and we obtain that the differential equation \eqref{eq16} is
equivalent to

\begin{equation}\label{eq17}\omega''(\xi)+\frac{\alpha}{2(n+m-1)}\omega'(\xi)\omega(\xi)^{\frac{4}{n+m-2}}=0.
\end{equation}
Integrating equation \eqref{eq17} we have

$$\omega'(\xi)+\frac{\alpha(n+m-2)}{2(n+m+2)(n+m-1)}\omega(\xi)^{\frac{n+m+2}{n+m-2}}=\beta,\hspace{0,3cm}\beta\in\mathbb{R}.
$$
Thus,
\begin{equation}-\int\frac{1}{\frac{\alpha(n+m-2)}{2(n+m+2)(n+m-1)}\omega(\xi)^{\frac{n+m+2}{n+m-2}}-\beta}d\omega=\xi+\nu\nonumber
\end{equation}
and then
\begin{equation}(n+m-1)(n+m+2)\int\frac{\varphi d\varphi}{\alpha-\frac{2\beta(n+m-1)(n+m+2)}{n+m-2}\varphi^{\frac{n+m}{2}+1}}=\xi+\nu\nonumber
\end{equation}
which is equation \eqref{eqc} of Theorem \ref{eq5}. Then we prove the
necessary condition. Now a direct calculation shows us the converse
implication. This concludes the proof of Theorem \ref{eq5}.
\end{myproof}

\begin{myproof}[\textbf{Proof of Theorem \ref{eq7}}]Since $h'=0$ and $\lambda_{F}=\rho=0$ we have by equation \eqref{eq1} and \eqref{eq2} of Theorem \ref{eq4} that

\begin{equation}(n-1)(2\varphi\varphi''-n(\varphi')^{2})-2\frac{m}{f}(\varphi^{2}f''-(n-2)\varphi\varphi'f')-\frac{m(m-1)}{f^2}\varphi^{2}(f')^{2}=0\nonumber
\end{equation}
which is equivalent to

\begin{equation}\left(\frac{f'}{f}-\frac{(n-2)}{(m+1)}\frac{\varphi'}{\varphi}\right)^{2}+\frac{2}{m+1}\left(\frac{f'}{f}-\frac{(n-2)}{(m+1)}\frac{\varphi'}{\varphi}\right)^{'}+\frac{(n+m-1)}{m(m+1)^{2}}\left(n(\frac{\varphi'}{\varphi})^{2}-2\frac{\varphi''}{\varphi}\right)=0.\nonumber
\end{equation}
Consider
$z=\frac{f'}{f}-\frac{(n-2)}{(m+1)}\frac{\varphi'}{\varphi}$, then

\begin{equation}\label{eq20}z^{2}+\frac{2}{m+1}z'+\frac{(n+m-1)}{m(m+1)^{2}}\left(n(\frac{\varphi'}{\varphi})^{2}-2\frac{\varphi''}{\varphi}\right)=0.
\end{equation}

Now, recall that the Ricatti differential equation is a differential
equation of the form

\begin{equation}\label{eq21}z(\xi)'=p(\xi)+q(\xi)z(\xi)+r(\xi)z(\xi)^{2}
\end{equation}
where $p$, $q$ and $r$ are smooth functions on $\mathbb{R}$, and by
Picard theorem we have that the solutions of \eqref{eq21} is given
by

$$z(\xi)=z_{p}(\xi)+\frac{e^{\int P(\xi)d\xi}d\xi}{-\int r(\xi)e^{\int P(\xi)d\xi}d\xi+c}$$
where $P(\xi)=q(\xi)+2z_{p}(\xi)r(\xi)$, $z_{p}(\xi)$ is a
particular solution of \eqref{eq21} and $c$ is a constant.

Observe that \eqref{eq20} is a Ricatti differential equation with
\begin{equation}q(\xi)=0,\hspace{0,2cm} r(\xi)=-\frac{m+1}{2}\hspace{0,2cm} \text{and}\hspace{0,2cm}
p(\xi)=-\frac{(n+m-1)}{2m(m+1)}\left(n(\frac{\varphi'}{\varphi})^{2}-2\frac{\varphi''}{\varphi}\right)
\end{equation}
Then we obtain

$$\frac{f'(\xi)}{f(\xi)}=\frac{(n-2)}{(m+1)}\frac{\varphi'(\xi)}{\varphi(\xi)}+z_{p}(\xi)+\frac{e^{-(m+1)\int z_{p}(\xi)d\xi}}{\frac{m+1}{2}\int e^{-(m+1)\int z_{p}(\xi)d\xi}d\xi+c}$$
and thus,
$$f=\varphi^{\frac{n-2}{m+1}}e^{\int z_{p}d\xi}\left(\int e^{-(m+1)\int z_{p}d\xi}d\xi+\frac{2}{m+1}C\right)^{\frac{2}{m+1}}$$
where $z_{p}(\xi)$ is a particular solution of \eqref{eq20}. This
expression is equation \eqref{eq22} of Theorem \ref{eq:02}.

Now, since $h'=0$, we have that $h(\xi)=constant$, which is equation
\eqref{eq23} of theorem 1.5. Then we prove the necessary condition.
Now a direct calculation shows us the converse implication. This
concludes the proof of Theorem \ref{eq:02}.
\end{myproof}

\begin{myproof}[\textbf{Proof of Theorem \ref{eq6}}]In this case, since $\lambda_{F}=\rho=0$, we have
by differential equation \eqref{eq:10} and \eqref{eq10} that

$$h(\xi)=\alpha\int\frac{1}{\varphi^{2}(\xi)}d\xi$$
for some $\alpha\neq0$ and $f$, $\varphi$ are arbitrary.

Then we prove the necessary condition. Now a direct calculation
shows us the converse implication. This concludes the proof of
Theorem \ref{eq6}.

\end{myproof}

\end{section}

\vskip0.8cm
\noindent
{Willian Isao Tokura} (e-mail: williamisaotokura@hotmail.com)\\[2pt]
Instituto de Matem\'atica e Estat\'istica\\
Universidade Federal de Goi\'as\\
74001-900-Goi\^ania-GO\\
Brazil\\

\noindent{Levi Adriano } (e-mail: levi@ufg.br)\\[2pt]
Instituto de Matem\'atica e Estat\'istica\\
Universidade Federal de Goi\'as\\
74001-900-Goi\^ania-GO\\
Brazil\\

\noindent{Romildo da Silva Pina} (e-mail: romildo@ufg.br)\\[2pt]
Instituto de Matem\'atica e Estat\'istica\\
Universidade Federal de Goi\'as\\
74001-900-Goi\^ania-GO\\
Brazil

\end{document}